\documentclass[a4paper,12pt]{article}
\usepackage{amssymb}
\usepackage{amsmath}
\usepackage{amsthm}
\usepackage{mathrsfs}
\usepackage[frame,curve,knot]{xy}
\xyoption{all}
\def\dj{d\kern-.30em\raise1.25ex\vbox{\hrule width .3em height .03em}}
\def\Dj{\rlap{\kern-.70em\raise0.75ex\vbox{\hrule width .3em height .03em}}} 

\newtheorem{theorem}{Theorem}[section]
\newtheorem{proposition}[theorem]{Proposition}
\newtheorem{lemma}[theorem]{Lemma}
\newtheorem{cor}[theorem]{Corollary}

\newtheorem{definition}[theorem]{Definition}
\newtheorem{example}[theorem]{Example}

\newtheorem{remark}[theorem]{Remark}

\newcommand{\bfrak}{\mathfrak{b}}
\newcommand{\cfrak}{\mathfrak{c}}
\newcommand{\dfrak}{\mathfrak{d}}

\newcommand{\chara }{\mbox{char}}

\newcommand{\cqg}{\mathbb{C}_q[G]}

\newcommand{\C }{\mathbb{C}}

\newcommand{\depth}{\mathrm{depth}}

\newcommand{\Ext}{\mathrm{Ext}}

\newcommand{\gfrak}{\mathfrak{g}}
\newcommand{\gkdim}{\mbox{GKdim}\,}

\newcommand{\hfrak}{\mathfrak{h}}

\newcommand{\hght}{\mathrm{ht}}

\newcommand{\Iqwld}{{\mathfrak I}_q(w,\lambda)}
\newcommand{\field}{\mathbb{K}}

\newcommand{\kqg}{\field_q[G]}

\newcommand{\N }{\mathbb{N}}

\newcommand{\oqgmn}{{\mathcal O}_q(G_{m,n}(\field))}
\newcommand{\oqmnm}{{\mathcal O}_q(M_{m,n}(\field))}

\newcommand{\ot }{\otimes }

\newcommand{\pfrak}{\mathfrak{p}}

\newcommand{\Ppluspi}{P^+(\pi)}

\newcommand{\Q}{\mathbb{Q}}

\newcommand{\rh }{\hat{R}}

\newcommand{\slfrak}{\mathfrak{sl}}
\newcommand{\sofrak}{\mathfrak{so}}

\newcommand{\sqgp}{S_q[G/P_S]}
\newcommand{\sqwgp}{S^w_q[G/P_S]}

\newcommand{\ubarm}{\overline{U}_-}

\newcommand{\uqg }{U_q(\mathfrak{g})}

\newcommand{\uqps }{U_q(\mathfrak{p}_S)}

\newcommand{\uqbp }{U_q(\mathfrak{b}^+)}

\newcommand{\uw }{\underline{w}}
\newcommand{\wght }{\mathrm{wt}}

\newcommand{\wurz }{\pi }

\newcommand{\vep }{\varepsilon }

\newcommand{\WSlen }{W^S_{n,\ge} }
\newcommand{\WSlenw }{W^S_{n,\ge,w} }

\newcommand{\WSlew }{W^S_{\le w} }

\newcommand{\Z}{\mathbb{Z}}

\title{The AS-Cohen-Macaulay property for quantum flag
  manifolds of minuscule weight}

\author{Stefan Kolb}

\date{July 10, 2007}

\begin{document}

\maketitle

\begin{abstract}
It is shown that quantum homogeneous coordinate rings of generalised flag manifolds
corresponding to minuscule weights, their Schubert
varieties, big cells, and determinantal varieties are
AS-Cohen-Macaulay. The main ingredient in the proof is the notion of a
quantum graded algebra with a straightening law, introduced by
T.H.~Lenagan and L.~Rigal [J.~Algebra {\bf 301} (2006), 670-702]. Using Stanley's
Theorem it is moreover shown that quantum generalised flag manifolds
of minuscule weight and their big cells are AS-Gorenstein. 
\end{abstract}

\noindent{\bf 2000 MSC:} 16E65, 16W35, 20G42

\noindent{\bf Keywords:} Quantum flag manifolds, straightening laws, Cohen-Macaulay, Gorenstein

\section{Introduction}
The theory of quantum groups provides noncommutative analogues
of coordinate rings of affine algebraic groups and their homogeneous
spaces. It is a recurring theme that these algebras resemble their
undeformed counterparts with respect to properties that allow a Lie
theoretic formulation. As pointed out in
\cite{a-GoodZhang07} this similarity is also a guiding principle for
homological properties of quantised coordinate rings. Notably
Cohen-Macaulay, (Auslander-) Gorenstein, and (Auslander-) regular
conditions should be reflected in the quantum setting. Indeed, for the generic
quantised coordinate ring $\cqg$ this is verified in
\cite[Theorem 0.1]{a-GoodZhang07}, \cite[Proposition 6.5]{a-BrownZhang}.

The homogeneous coordinate rings of generalised flag manifolds allow
straightforward quantum analogues as $\N_0$-graded algebras
\cite{a-LaksResh92}, \cite{a-Soib92}. A
graded analogue of the Cohen-Macaulay property is the AS-Cohen
Macaulay property introduced in  \cite[p.~674]{a-vdB97}, \cite{a-Jorg99}.
For quantum generalised
flag manifolds it is crucial to establish the AS-Cohen-Macaulay
property for the following reason. The commutation
relations of quantum algebras are governed by the universal
$R$-matrix. The construction of the $R$-matrix implies that quantum
generalised flag manifolds have enough normal elements in the sense of
\cite{a-Zhang97}. Hence the methods developed in that paper may be
applied. In particular, given AS-Cohen-Macaulayness, one may use
Stanley's Theorem \cite[Theorem 6.2]{a-JorgZhang00} to verify
AS-Gorensteinness, which in turn implies the Auslander-Gorenstein and
the Cohen-Macaulay property \cite[Theorem~0.2]{a-Zhang97}. Thus
establishing AS-Cohen-Macaulayness is an essential step towards other
desirable homological properties.

In the present paper the AS-Cohen-Macaulay property is verified for
quantum homogeneous coordinate rings of generalised flag manifolds
corresponding to minuscule weights, their Schubert
varieties, big cells, and determinantal varieties. Using Stanley's
Theorem it is moreover shown that quantum generalised flag manifolds
of minuscule weight and their big cells are AS-Gorenstein. 

The main ingredient in the proof is the notion of a quantum graded
algebra with a straightening law (ASL) introduced in \cite{a-LenRi06},
in generalisation of a notion from commutative algebra developed by
C.~De Concini, D.~Eisenbud, and C.~Procesi
(e.g.~\cite{a-DEP82}, \cite{b-BV88}). In their paper T.H.~Lenagan
and L.~Rigal prove that any quantum
graded ASL on a wonderful poset is AS-Cohen-Macaulay. They
apply this result to quantum Grassmannians, their Schubert
varieties, quantum matrices and the corresponding quantum
determinantal varieties. Moreover, in \cite{a-LenRi06} AS-Gorensteinness
is established via Stanley's Theorem. 

The original motivation for the present work was to understand the
arguments of \cite{a-LenRi06} in a more representation theoretic
setting. In particular a consequent use of the Bruhat order, standard
monomial theory, and the properties of
the universal $R$-matrix shorten the arguments and allow
generalisation to all flag manifolds corresponding to minuscule weights.
One slight remark of caution, however, is in place. Following
\cite{b-Jantzen96} we work over an arbitrary field $\field$ and with a deformation
parameter $q\in \field\setminus \{0\}$ which is not a root of
unity. No attempt is made to include the case where $q$ is a root of
unity, although the results of \cite{a-LenRi06} also hold in this
setting.

The paper is organised as follows. In Section \ref{prelim} we fix the
usual notations for semisimple Lie algebras and recall the notions of
wonderful posets, minuscule weights, and Seshadri's classical result
on the standard monomial theory for flag manifolds of minuscule
weight. In Section \ref{QGaQFM} we recall notions and results from the
theory of quantum groups and the definition and properties of quantum
homogeneous coordinate rings of flag manifolds, their Schubert
varieties, big cells, and determinantal varieties. We also consider
standard monomial theory in the minuscule case in Section
\ref{standardMonMin}.
Generically these results are well known. Following \cite{b-Jantzen96} however, we
take great care that all our statements hold over an arbitrary field $\field$
and for any $q\in \field\setminus\{0\}$ which is not a root of
unity. Moreover, for $\gfrak=\slfrak_n$ we identify the representation
theoretic definition of quantum homogeneous coordinate rings of Grassmann
manifolds with the definition in terms of quantum matrices and their quantum
minors. 

With these preliminary considerations and \cite[Theorem
  2.2.5]{a-LenRi06} at hand we give a short proof of the
AS-Cohen-Macaulay property in Section \ref{ASCM-property}. Moreover,
the AS-Gorenstein property for quantum homogeneous coordinate rings
of flag manifolds of minuscule weight and their big cells is
established in \ref{ASG}. This argument does not involve any quantum
considerations. 

Finally, in \ref{example} the AS-Cohen Macaulay and AS-Gorenstein properties
are established for the simplest family of quantum flag manifolds
corresponding to non-minuscule fundamental weights.

\section{Preliminaries}\label{prelim}
Let $\C$, $\Q$, $\Z$, $\N$, and $\N_0$ denote the complex numbers, the
rational numbers, the integers, the positive integers, and the nonnegative integers, respectively.

\subsection{Semisimple Lie algebras}
 Let $\gfrak$ be
a finite dimensional complex semisimple Lie algebra of rank $r$ and
$\hfrak\subseteq\gfrak$ a fixed Cartan subalgebra. Let
$\Delta\subseteq \hfrak^\ast$  denote the root system  associated with
$(\gfrak,\hfrak)$. Choose an ordered basis
$\wurz=\{\alpha_1,\dots,\alpha_r\}$ of simple roots for $\Delta$ and
let $\Delta^+$ (resp.~$\Delta^-$) be the set of positive 
(resp.~negative) roots with respect to $\wurz$. Identify $\hfrak$ with
its dual via the Killing form. The induced
nondegenerate symmetric bilinear form on $\hfrak^*$ is denoted by
$(\cdot,\cdot)$.
For $\alpha\in \Delta$ let $\alpha^\vee=2\alpha/(\alpha,\alpha)$
denote the corresponding coroot.
Let $\omega_i\in\hfrak^\ast$,
$i=1,\dots,r$, be the fundamental weights with respect to $\pi$ and let 
$P^+(\wurz)$ denote the set of dominant integral weights.
Moreover, let $\leq$ denote the standard partial ordering on $\hfrak^*$. In
particular, $\mu\leq\gamma$ if and only  if $\gamma-\mu\in \N_0\pi$.
For $\mu\in P^+(\pi)$ let $V(\mu)$ denote the finite dimensional 
irreducible $\gfrak$-module of highest weight $\mu$. We will write
$\Omega(V)$ to denote the set of weights of a finite dimensional $\gfrak$-module $V$.

Let $G$ denote the connected, simply connected complex  Lie group with Lie algebra $\gfrak$.
For any subset $S\subset \pi$ let $P_S\subset G$ be the
corresponding standard parabolic subgroup with Lie algebra
$$\pfrak_S=\hfrak\oplus\bigoplus_{\alpha\in \Delta^+\cup (\Z S
  \cap\Delta^-)} \gfrak_\alpha,$$  
where $\gfrak_\alpha$ denotes the root space corresponding to $\alpha\in \Delta$.

Let $W$ be the Weyl group associated to the root system $\Delta$. 
For any $\alpha\in \Delta^+$ let $s_\alpha\in W$ denote the reflection
on the hyperplane orthogonal to $\alpha$ with respect to $(\cdot,\cdot)$. 
Moreover, for any subset $S\subset \pi$ we write $W_S\subset W$ to denote the subgroup
generated by the reflections corresponding to simple roots in $S$, and
we write $W^S$ for the set of minimal length coset representatives  of $W/W_S$ in $W$.

For our purposes it is essential that the Weyl group is a poset with
the Bruhat order. For $w,w'\in W$ write $w\rightarrow w'$ if there exists
$\alpha\in \Delta^+$ such that $w=s_\alpha w'$ and $l(w')=l(w)+1$
where $l$ denotes the length function on $W$. The Bruhat
order $\le $ on $W$ is then given by the relation
\begin{align*}
  w\le w' \Leftrightarrow \
  & \text{there exists } n\ge 1,\ w_2,\ldots ,w_{n-1}\in W,\\
  & \text{such that }w=w_1\rightarrow w_2\rightarrow \ldots
  \rightarrow w_n=w'.
\end{align*}
Note that we are using the same symbol $\le$ for the standard partial
ordering on the weight lattice and for the Bruhat order. We hope that
this does not lead to confusion.

\subsection{Wonderful posets}
Let $(P,\le)$ be a poset. Inspired by the notation for the Bruhat
order we write  $x\rightarrow y$  for two elements $x,y\in P$ if $x<y$
and there does not exist an element $z\in P$ such that $x<z<y$. In
this case $y$ is called an upper neighbour of $x$. For our purposes
the following slightly simplified notion of a definition from \cite[Chapter 5.D]{b-BV88} is sufficient.
\begin{definition}
 A poset $(P,\le)$ with a smallest and a greatest element is called
  wonderful if for any $x,y,z,u\in P$ such that  $z \rightarrow x\le u$ and $z \rightarrow
  y\le u$ there exists an element  $w \in P$ such that $x \rightarrow w$, $y
  \rightarrow w$, and $w\le u$.
\end{definition} 
In the examples we are interested in $P$ is always a subset of the Weyl
group and $\le$ is the Bruhat order. These posets will have
a smallest and a greatest element.
Note also that the Weyl group itself
with the Bruhat order is in general not wonderful.

\begin{example}\upshape Consider the symmetric group $S_4$ with standard generators
$s_1=(12)$, $s_2=(23)$, $s_3=(34)$. Then the elements
$x=s_1s_2s_3$, $y=s_3s_2s_1$, $z=s_1s_3$ satisfy the relations $z
\rightarrow x$ and $z \rightarrow y$. However, one easily verifies
that there exists no element $w$ such that
$x \rightarrow w$ and $y \rightarrow w$. Hence $S_4$ is not wonderful
with the Bruhat order.
\end{example}

\subsection{Minuscule weights}
  Let $\lambda\in \Ppluspi\setminus\{0\}$ be a minuscule weight in the sense of
\cite[page 72]{b-Humphreys}, i.e.~ $(\lambda,\alpha^\vee)\in
\{0,1\}$ for all positive roots $\alpha$. Humphreys calls such weights
 minimal. Recall that this condition
implies that the set of weights $\Omega(V(\lambda))$ consists of one
single orbit under the action of the Weyl group. In particular all
weight spaces of $V(\lambda)$ are one-dimensional. Moreover,
$\lambda=\omega_s$ for some simple root $\alpha_s\in \pi$. A
complete list of the possible $\omega_s$ can be found also in
\cite[page 72]{b-Humphreys}.
\begin{proposition}\label{wonderprop}
   Let $\lambda=\omega_s$ be a minuscule weight and 
  $S=\pi\setminus \{\alpha_s\}$. Then the following hold:
  \begin{enumerate}
  \item Assume that $w,w'\in W^S$ satisfy $w\rightarrow w'$. Then
  $w'=s_i w$ for some $\alpha_i\in \pi$.
  \item Assume that $w,s_iw,s_jw \in W^S$ and
  $l(s_iw)=l(s_jw)=l(w)+1$ for some $\alpha_i,\alpha_j\in \pi$. Then
  $s_is_jw=s_js_i w\in W^S$ holds.
  If moreover $\alpha_i\neq\alpha_j$ then $l(s_is_jw)=l(w)+2$.
  \item Assume that $w,s_iw,s_jw \in W^S$ and
  $l(s_iw)=l(s_jw)=l(w)-1$ for some $\alpha_i,\alpha_j\in \pi$. Then
  $s_is_jw=s_js_i w\in W^S$ holds.
  If moreover $\alpha_i\neq\alpha_j$ then $l(s_is_jw)=l(w)-2$.
  \item  Two elements  $w,w'\in W^S$ satisfy $w\le w'$
  if and only if $w\lambda\ge w'\lambda$.
  \end{enumerate}
\end{proposition}

\noindent {\bf Proof:} We leave it to the reader to verify these elementary
  facts. One way around the proof is to construct the representation
  $V(\omega_s)$ explicitly, as it is done for instance
  in \cite[5A.1]{b-Jantzen96} in the quantum case.
  $\blacksquare$

\medskip

 For any $S\subseteq \pi$ define $W^S_{\le w}=\{w'\in W^S\,|\, w'\le w\}$. 
The first, the second, and the last statement of the above proposition
immediately imply the following result.
\begin{cor}\label{wonderful} 
  Let $\lambda=\omega_s$ be a minuscule weight,
  $S=\pi\setminus \{\alpha_s\}$,
  and $w\in W^S$. Then the poset $(W^S_{\le w}, \le)$ is wonderful.
  In particular the poset $(W^S,\le)$ is  wonderful.
\end{cor}

\begin{remark}\label{invBruRemark}\upshape
  In our later conventions it will be advantageous to consider
  $W^S_{\le w}$ with the inverse Bruhat order $\ge$. The first,
  the third, and the last statement of the above proposition imply that for minuscule
  $\omega_s$ the poset $(W^S_{\le w},\ge)$ is also wonderful. $\square$
\end{remark}

\subsection{A result from standard monomial theory}\label{standard}
Among the many useful applications of standard monomial theory are
character formulae for irreducible
representations $V(\lambda)$, $\lambda\in \Ppluspi$, in terms of the
Weyl group and its subsets $W^S$ where $S\subset \pi$. 
We introduce the following notation to describe characters for
representations of highest weight $n\omega_s$ where
$n\in \N_0$ and $\omega_s$ is minuscule. For any subset $S$ of $\pi$
and $n\in \N_0$ define
\begin{align*}
  W^S_{n,\ge}:=\{(w_1,\dots,w_n)\in (W^S)^n\,|\, w_{j}\ge w_{j+1}
  \mbox{ for } j=1,\dots,n-1\}.
\end{align*}
For any weight $\mu$ let $e^\mu$ denote the corresponding element in
the group algebra of the weight lattice.
Moreover, we write $\mbox{char}(V)$ to denote the character of a finite dimensional
$\gfrak$-module $V$.
The next result follows from the standard monomial theory
developed in \cite{a-seshadri78}. 
\begin{proposition}\label{charProp}
  Let $\lambda=\omega_s$ be a minuscule weight,
 $S=\pi\setminus\{\alpha_s\}$, and $n\in \N_0$.
 Then 
  \begin{align*}
    \mbox{\upshape char}(V(n\lambda))=\sum_{(w_1, \cdots,w_n)\in W^S_{n,\ge}} 
    e^{w_1 \lambda}e^{w_2\lambda} \cdots e^{w_n\lambda}.
  \end{align*}
\end{proposition}

\section{Quantum groups and quantum flag manifolds}\label{QGaQFM}
  In this section we recall well known objects from the theory of quantum
  groups, in particular quantised algebras of functions on flag
  manifolds and Schubert varieties. All properties used to prove
  the main result hold over an arbitrary field $\field$ and for any
  deformation parameter $q\in\field\setminus \{0\}$ which is
  not a root of unity.

\subsection{Quantum groups}
We recall some well know definitions and facts
along the lines of \cite{b-Jantzen96}. Let $\field$ be an arbitrary
field and let $q\in \field\setminus\{0\}$ be not a root of
unity. Choose $N\in \N$ minimal such that $(\omega_i,\omega_j)\in
\Z/N$ and assume that $\field$ also contains $q^{1/N}$.
The quantum universal enveloping algebra $\uqg$ is the $\field$-algebra
generated by elements  $E_i, F_i, K_i, K_i^{-1}$, $i=1,\dots,r$, and relations given
for instance in \cite[4.3]{b-Jantzen96}.
Let $\uqbp$ denote the subalgebra generated by $E_i, K_i, K_i^{-1}$
for $i=1,\dots,r$.
For $\lambda\in \Ppluspi$ we write $V(\lambda)$ to denote the
irreducible $\uqg$-module of highest weight $\lambda$. In particular
$V(\lambda)$ contains a highest weight vector $v_\lambda$ which satisfies
\begin{align*}
  E_i v_\lambda=0, \qquad K_i v_\lambda= q^{(\lambda,\alpha_i^\vee)}v_\lambda.
\end{align*}
We hope it will be clear from the context whether
$V(\lambda)$ denotes a $\gfrak$-module or an $\uqg$-module. The
character of the $\gfrak$-module $V(\lambda)$ coincides with the
character of the corresponding $\uqg$-module. Hence we may apply
Proposition \ref{charProp} and the Weyl character formula also in the quantum case.

Note that the dual $\field$-vector space $V(\lambda)^\ast$ is a right
$\uqg$-module with action given by $fu(v)=f(uv)$ for all $f\in
V(\lambda)^\ast, v\in V(\lambda), u\in \uqg$.
For any $v\in V(\lambda)$ and $f\in V(\lambda)^\ast$ we write
$c^\lambda_{f,v}$ to denote the linear functional on $\uqg$ defined by
$c^\lambda_{f,v}(u)=f(uv)$ for all $u\in \uqg$. The quantum coordinate
ring $\kqg$ is defined to be the linear span of the all matrix
coefficients $c^\lambda_{f,v}$ for all $\lambda\in \Ppluspi$. It is a
$\uqg$-bimodule algebra with right and left action given by
$u c^\lambda_{f,v} u'=c^\lambda_{fu',uv}$ for all
$u,u'\in \uqg$.

It is an important observation due to V.~Drinfeld that the category of
type one representations of $\uqg$, i.e.~of finite
direct sums of representations $V(\lambda)$, is a braided monoidal
category \cite[7.3--7.8]{b-Jantzen96}. Moreover, for any two type one representations $V,W$ and any
weight vectors $v\in V$, $w\in W$ the braiding
$\rh_{V,W}:V\ot W \rightarrow W\ot V$ satisfies
\begin{align}\label{R-mat-prop}
  \rh_{V,W}(v\ot w)=q^{(\wght(v),\wght(w))}w\ot v +\sum_i w_i \ot v_i
\end{align}
where $w_i,v_i$ are weight vectors satisfying $\wght(w_i){<}\wght(w)$
and $\wght(v_i){>}\wght(v)$. This property follows from the
construction of $\rh_{V,W}$.
Formula (\ref{R-mat-prop}) implies
that for any $\lambda\in \Ppluspi$ the vector $v_\lambda \ot
v_\lambda\in V(\lambda)\ot V(\lambda)$ is an eigenvector of 
$\rh_{V(\lambda),V(\lambda)}$ with eigenvalue
$q^{(\lambda,\lambda)}$.

\begin{remark}\label{LR-remark} \upshape
  It is possible to translate \cite[Lemma 2.16]{a-LaksResh92} into the
  language of \cite[Chapter 7]{b-Jantzen96}. This implies, also for
  general $\field$ and $q$ not a root of unity, that  
  $V(2\lambda)\subset  V(\lambda)\ot V(\lambda)$ is the full
  eigenspace of $\rh_{V(\lambda),V(\lambda)}$ corresponding to the eigenvalue
  $q^{(\lambda,\lambda)}$. $\square$
\end{remark}

\begin{remark} \upshape
  Note that \cite[Chapter 7]{b-Jantzen96} is formulated for
  $\chara(\field)=0$ and any $q$ which is transcendental over $\Q$. As pointed out in
  the introduction to that chapter, however, everything can be extended to the
  case of general characteristic and $q$ not a root of unity. To this
  end it is essential to note that the pairing
  \cite[6.12]{b-Jantzen96} remains nondegenerate. This follows from
  \cite[6.23, 8.30]{b-Jantzen96}. $\square$
\end{remark}

\subsection{Quantum flag manifolds}\label{q-flags}
For any $\lambda\in \Ppluspi$ define $S\subset \pi$ by
$S=S(\lambda)=\{\alpha_i \in \pi \,|\, (\lambda,\alpha_i)\neq 0\}$.  The
generalised flag manifold $G/P_S$ is a projective algebraic variety. It
can be embedded into projective space ${\mathbb P}(V(\lambda))$ by $G\ni
g\mapsto [gv_\lambda]$ where $v_\lambda\in V(\lambda)$ is the
highest weight vector and $[v]$ denotes the line represented by $v\in V(\lambda)$.
The homogeneous coordinate ring $S[G/P_S]$ with respect to this
embedding is isomorphic to $\bigoplus_{n=0}^\infty V(n\lambda)^\ast$
endowed with the Cartan multiplication
\begin{align}\label{cartan}
  V(n_1\lambda)^\ast \ot V(n_2\lambda)^\ast
  \rightarrow V((n_1+n_2)\lambda)^\ast.
\end{align}
These structures can immediately be translated into the quantum
setting. Indeed, following \cite{a-LaksResh92}, \cite{a-Soib92}, \cite{b-CP94} we define 
the quantised homogeneous coordinate ring $\sqgp$ of the generalised flag
manifold $G/P_S$ to be the subalgebra of $\kqg$ generated by the
matrix coefficients $\{c_{f,v_\lambda}^\lambda\,|\, f\in V(\lambda)^\ast\}$,
where $v_\lambda$ is a highest weight vector of the $\uqg$-module $V(\lambda)$.
We will freely identify $f\in V(\lambda)^\ast$ with the generator $c^\lambda_{f,v_\lambda}$.
It is immediate from the definition of the multiplication on $\kqg$
that the right $\uqg$-module algebra $\sqgp$ is again isomorphic to 
$\bigoplus_{n=0}^\infty V(n\lambda)^\ast$ endowed with the Cartan
multiplication (\ref{cartan}). Note that $\sqgp$ depends on $\lambda$
and not only on $S$ even if this is not made explicit in our notation.

Let $\rh^\ast_{V(\lambda),V(\lambda)}: V(\lambda)^\ast\ot
V(\lambda)^\ast \rightarrow V(\lambda)^\ast\ot
V(\lambda)^\ast $ denote the map dual to $\rh_{V(\lambda),V(\lambda)}$.
Recall from the previous subsection that the subspace $V(2\lambda)\subset V(\lambda)\ot V(\lambda)$ is
the eigenspace of $\rh_{V(\lambda),V(\lambda)}$ with corresponding eigenvalue
$q^{(\lambda,\lambda)}$. It follows that the generators $f,g\in
V(\lambda)^\ast $ of $\sqgp$ satisfy the relations
\begin{align}\label{flag-rels}
  m(\rh^\ast_{V(\lambda),V(\lambda)}(f\ot g))=q^{(\lambda,\lambda)}f g
\end{align}
where $m$ denotes the multiplication map.

\begin{remark}\label{gen-rels}\upshape
  It is moreover known \cite{a-Brav94} that $\sqgp$ is a quadratic
  algebra. Check carefully that in view of Remark \ref{LR-remark}
  Braverman's result also holds in general characteristic and for $q$
  not a root of unity.
  Hence the algebra
  $\sqgp$ is given in terms of generators and  relations by (\ref{flag-rels}).
  This is an instance of a quantum effect. In the commutative case the
  $R$-matrix becomes the twist of tensor factors, and there is no
  similar canonical way to encode the defining Pl\"ucker relations of
  $S[G/P_S]$. $\square$
\end{remark}
\begin{remark}\upshape
  To facilitate later reference we explicitly translate property (\ref{R-mat-prop}) into
  properties of $\rh^\ast_{V(\lambda),V(\lambda)}$. For $f\in V(\lambda)^\ast\setminus\{0\}$ write
  $\wght(f)=\mu$ if $f(v)=0$ for all weight vectors $v\in V(\lambda)$
  with $\wght(v)\neq \mu$. For $f,g\in V(\lambda)^\ast$ formula
  (\ref{R-mat-prop}) implies
  \begin{align}\label{R-mat-prop-dual}
    \rh_{V(\lambda),V(\lambda)}^\ast(f\ot g)=q^{(\wght(f),\wght(g))}g\ot f +\sum_i g_i \ot f_i
  \end{align}
  where $g_i,f_i\in V(\lambda)^\ast$ are weight vectors such that $\wght(g_i)<\wght(g)$
  and $\wght(f_i)>\wght(f)$. $\square$
\end{remark}

\begin{remark}\upshape
  Note that for $\gfrak=\slfrak_n$ and $\lambda=\omega_m$
  the ring $\sqgp$ coincides with the algebra $\oqgmn$
  considered in \cite{a-LenRi06}. Indeed, recall that by definition
  $\oqgmn$ is the subalgebra generated by quantum $(m\times m)$-minors
  in the algebra $\oqmnm$ of quantum $(m\times n)$-matrices. Here we
  consider $\oqmnm$ as a left $U_q(\slfrak_m)$-module algebra and a right
  $U_q(\slfrak_n)$-module algebra generated by $V(\omega_1)\ot
  V(\omega_1)^\ast$, where
  $V(\omega_1)$ and $V(\omega_1)^\ast$ are considered as left
  $U_q(\slfrak_m)$-module and right $U_q(\slfrak_n)$-module, respectively. With
  respect to the right $U_q(\slfrak_n)$-module algebra structure
  the space spanned by the quantum $(m\times m)$-minors is
  isomorphic to $V(\omega_m)^\ast$. The commutation relations of
  $\oqmnm$ are given in terms of $\rh$-matrices and induce the relations
  (\ref{flag-rels}) on quantum $(m\times m)$-minors. This follows from
  the naturality of the braiding and the hexagon equation. As $\oqmnm$ is a
  domain it follows from Remark \ref{gen-rels} that $\oqgmn$ is
  isomorphic to $\sqgp$ as a $U_q(\slfrak_n)$-module algebra.
  $\square$
\end{remark}

\begin{lemma}\label{noetherlem}
  For any $\lambda\in \Ppluspi$ and $S=S(\lambda)$ the algebra $\sqgp$
  is noetherian.
\end{lemma}
\noindent{\bf Proof:}
  Let $x_1,\dots,x_{\dim V(\lambda)}$ be a weight basis of
  $V(\lambda)^\ast$ such that $\wght(x_i)<\wght(x_j)$ implies
  $i<j$. It follows from the commutation relation
  (\ref{flag-rels}) and (\ref{R-mat-prop-dual}) that for all $1\le i\le j\le \dim
  V(\lambda)$ there exist scalars $q_{ij}\in \field\setminus \{0\}$ and
  $\alpha_{ij}^{st}\in \field$ such that
  \begin{align*} 
    x_ix_j=q_{ij} x_jx_i +\sum_{s=1}^{j-1} 
    \sum_{t=1}^{\dim V(\lambda)}\alpha_{ij}^{st} x_s x_t.
  \end{align*}
  Now the lemma follows from \cite[Proposition I.8.17]{b-BG02}.
$\blacksquare$

\subsection{Quantum Schubert Varieties}\label{schubert}

For any $\mu\in \Ppluspi$ and $w\in W$ let  $v_{w\mu}\in V(\mu)$ be a weight vector of weight
$w\mu$ and let $V_w(\mu)=\uqbp v_{w\mu}$ denote the corresponding
Demazure module. We write
$V_w(\mu)^\perp$ to denote the orthogonal complement of $V_w(\mu)$ in
$V(\mu)^\ast$. 

Quantum analogues of Schubert varieties were
introduced and studied in \cite{a-LaksResh92} (cp.~also
\cite{a-Soib92}, \cite{b-CP94}). Let $S$ and $\lambda$ be as in
Subsection \ref{q-flags}, and $w\in W$.
Let $\Iqwld$ be the two sided ideal in
$\sqgp$ generated by all matrix elements $c^\lambda_{f,v_\lambda}$
such that $f\in V_w(\lambda)^\perp$. We call the quotient algebra
$S_q^w[G/P_S]=\sqgp/\Iqwld$ the quantised algebra of functions on the
Schubert variety corresponding to $w$. Note that
$\Iqwld$ is right $\uqbp$-invariant and hence $S_q^w[G/P_S]$ is a right
$\uqbp$-module algebra. 

Following \cite[10.1.8]{b-Joseph} define
\begin{align*}
 Q_w=\sum_{n\in \N_0}V_w(n\lambda)^\perp. 
\end{align*}
We write $\pi_w:\sqgp\rightarrow \sqgp/Q_w$ to denote the canonical
projection. By \cite[10.1.8]{b-Joseph} the subspace $Q_w\subset \sqgp$
is a completely prime ideal invariant under the right action of
$\uqbp$. Note that the proof of this result given in \cite{b-Joseph} is also valid
for an arbitrary field $\field$ an any $q\in \field\setminus\{0\}$ which is not a
root of unity. It is moreover proved in \cite{a-Joseph95} that the ideals $\Iqwld$
and $Q_w$ coincide for transcendental $q$. In Subsection \ref{standardMonMin}
we give an elementary proof of this fact for minuscule $\lambda$ and
all $q$ which are not a root of unity.

\subsection{Quantum big cells and determinantal varieties}
Let $f_e\in V(\lambda)^\ast$ denote the up to a scalar factor uniquely
determined element of weight $\lambda$. The commutation relations
(\ref{flag-rels}) and (\ref{R-mat-prop-dual}) imply that the element
$f_e$ commutes up to a power of $q$ with all generators of $\sqgp$.
Hence we can form the localisations $\sqgp(f_e^{-1})$ and
$\sqwgp(f_e^{-1})$ of the quantum coordinate rings of flag manifolds and
Schubert varieties with respect to the multiplicative set
$\{f_e^n\,|\, n\in \N\}$. As $f_e$ is homogeneous of degree one these
localisations are $\Z$-graded and we let $\sqgp(f_e^{-1})_0$ and
$\sqwgp(f_e^{-1})_0$ denote the homogeneous components of degree
zero. We call $\sqgp(f_e^{-1})_0$ the quantum big cell of
$G/P_S$. Moreover, we call $\sqwgp(f_e^{-1})_0$ for $w\in
W$ quantum determinantal varieties. 

The algebra $\sqgp(f_e^{-1})_0$ can be realised as the graded dual of a
suitable coalgebra. This observation is not necessary for the proof of
Corollary \ref{minuscule-ASCM} but we include it here as we believe it to be
noteworthy in itself. Moreover, the construction explains how we want to consider
$\sqgp(f_e^{-1})_0$ and $\sqwgp(f_e^{-1})_0$ as graded algebras. It is
also used to establish the AS-Gorenstein property for
$\sqgp(f_e^{-1})_0$ in Proposition \ref{cell-ASGor}.

We introduce some notations to formulate the desired
duality. First define a relative height function on the root lattice by
\begin{align*}
  \hght_S:\Z\pi\rightarrow \Z,\qquad \hght_S\Big(\sum_{\alpha_i\in
  \pi}n_i\alpha_i\Big)=\sum_{\alpha_i\notin S}n_i.
\end{align*}
Note that the algebra $\sqgp(f_e^{-1})_0$ is $\N_0$-graded if one
defines $\deg(f_e^{-1}f_\mu)=\hght_S(\lambda-\mu)$ for any generator
$f_\mu\in V(\lambda)^\ast$ of weight $\mu$.
Moreover, define $\uqps$ to be the subalgebra of $\uqg$ generated by
the elements
$\{K_i^{\pm 1},E_i, F_j,\,|\,i=1,\dots,r\mbox{ and } \alpha_j\in S\}$ and let
$\uqps^+=\uqps\cap \ker(\vep)$ denote the augmentation ideal. Then 
\begin{align*}
  \ubarm:=\uqg/\uqg \uqps^+
\end{align*}
is a $\N_0$-graded coalgebra, where the grading is given by
$\deg(u)=\hght_S(\mu)$ if $u\in (\ubarm)_{-\mu}$.

Note that both $\sqgp(f_e^{-1})_0$ and $\ubarm$ have finite
dimensional homogeneous components. Let $\delta(\ubarm)$ denote the
graded dual algebra of $\ubarm$. The left $\uqg$-module structure on
$\ubarm$ induces the structure of a right $\uqg$-module algebra on
$\delta(\ubarm)$. 
\begin{proposition}\label{dualprop}
  There is an isomorphism of right $\uqg$-module algebras $\sqgp(f_e^{-1})_0\cong\delta(\ubarm)$.
\end{proposition}
\noindent{\bf Proof:}
Inside $\sqgp(f_e^{-1})$ one has an inclusion
\begin{align*}
  f_e^{-n}V(n\lambda)^\ast\subseteq  f_e^{-(n+1)}V((n+1)\lambda)^\ast
\end{align*}
and hence $\sqgp(f_e^{-1})_0$ can be written as a direct limit
\begin{align}\label{directLim}
 \sqgp(f_e^{-1})_0\cong \lim_{n\rightarrow \infty} f_e^{-n}V(n\lambda)^\ast.
\end{align}
Moreover, the canonical pairing $\kqg\ot \uqg\rightarrow \field$
induces a pairing 
\begin{align}\label{pair1}
  \sqgp(f_e^{-1})_0\ot \uqg \rightarrow \field.
\end{align}
It follows from (\ref{directLim}) that (\ref{pair1}) induces a
nondegenerate pairing 
\begin{align}\label{pairnondeg}
  \sqgp(f_e^{-1})_0\ot \ubarm \rightarrow \field.
\end{align}
between the algebra $ \sqgp(f_e^{-1})_0$ and the coalgebra $\ubarm$ which respects the
$\N_0$-grading of both $ \sqgp(f_e^{-1})_0$ and $\ubarm$. Therefore the
algebra $ \sqgp(f_e^{-1})_0$ coincides with the graded dual
$\delta(\ubarm)$.
$\blacksquare$

\begin{remark} \upshape
  For $\gfrak=\slfrak_n$ and $\lambda=\omega_s$ it is shown in \cite[Corollaries 1 and 2]{a-heko04} that
  the graded algebra $\delta(\ubarm)$ coincides with the graded algebra ${\mathcal
  O}_q(M_{s,{n-s}}(\field))$ of quantum $(s\times(n-s))$-matrices. In
  view of the above proposition there is also an isomorphism of graded
  algebras between $\sqgp(f_e^{-1})_0$ and ${\mathcal
  O}_q(M_{s,{n-s}}(\field))$. Analogous results hold for all
  cominuscule weights. Note that a
  direct proof of the isomorphism  $\sqgp(f_e^{-1})_0\cong{\mathcal
  O}_q(M_{s,{n-s}}(\field)) $ is also given in \cite{a-KLR04}.

  For $\gfrak=\slfrak_n$ and $\lambda=\omega_s$ one can moreover
  verify, that for suitable $w_t\in W^S$ the algebra
  $S^{w_t}[G/P_S](f_e^{-1})_0$ coincides with the algebra ${\mathcal
  O}_q(M_{s,{n-s}}(\field))/{\mathcal I}_t$ defined for instance in
  \cite[3.5]{a-LenRi06}. This explains the name quantum determinantal variety.
\end{remark}

\subsection{Standard monomials for minuscule weights}\label{standardMonMin}

Recall the notation $\WSlen$ introduced in Section \ref{standard}.
Moreover, for $n\in \N_0$ and $w\in W^S$ define
\begin{align*}
  \WSlenw&=\{(w_1,w_2,\dots,w_n)\in \WSlen\,|\,w_1\le w\}.
\end{align*}
For $\uw=(w_1,\dots,w_n)\in \WSlen$ we define standard monomials
$f_{\uw}\in V(n\lambda)^\ast \subset \sqgp$ by $f_{\uw}=f_{w_1}\dots
f_{w_n}$. Here $f_w\in V(\lambda)^\ast$ denotes the up to a scalar
factor uniquely determined element of weight $w\lambda$.

\begin{proposition} \label{standardMinProp}
  Let $\lambda=\omega_s$ be minuscule, $S=\pi\setminus \{\alpha_s\}$, and $w\in W^S$.
  \begin{enumerate}
    \item The equality $\Iqwld=Q_w$ holds, i.e.~ $\pi_w(\sqgp)=S_q^w[G/P_S]$.
    \item The standard monomials $\{f_{\uw}\,|\,\uw\in \WSlen\}$ form
       a basis of $V(n\lambda)^\ast\subset \sqgp$.
    \item The standard monomials $\{f_{\uw}\,|\,\uw\in \WSlenw\}$ form
      a basis of $\pi_w(V(n\lambda)^\ast)\subset \pi_w(\sqgp)$.
  \end{enumerate}
\end{proposition}

\noindent{\bf Proof:} Statement 1.~is equivalent to
$\Iqwld\cap V(n\lambda)^\ast=V_w(n\lambda)^\perp$ for all $n\in
\N_0$. Note that all statements of the above proposition hold for
$n=0,1$. We proceed by induction on $n$.  

We first prove linear independence of the standard monomials
$f_{\uw}\in \sqgp$ for $\uw\in \WSlen$. Assume there exists a
nontrivial linear combination
\begin{align*}
  \sum_{\uw\in \WSlen} c_{\uw} f_{\uw}=0.
\end{align*}
Choose $w\in W^S$ minimal such that $c_{\uw}\neq 0$ for some
$\uw=(w,w_2\dots,w_n)$. Then in $\sqgp/Q_{w}$ one obtains a relation
\begin{align}\label{nullteiler}
  f_{w}\left(\sum_{\uw'\in W^S_{n-1,\ge, w}} c_{\uw'}f_{\uw'}\right)=0
\end{align}
where not all coefficients $c_{\uw'}$ vanish. As $\sqgp/Q_{w}$ is an
integral domain one obtains a contradiction to 3.~for $n{-}1$.
Hence the standard monomials of length $n$ are linearly independent. By
Proposition \ref{charProp} they generate a space of dimension
$\dim(V(n\lambda))$ and hence form a basis. This proves 2.

The third claim is proved in a very similar fashion, but now one
also performs induction over $l(w)$. Assume there exists a
nontrivial linear combination
\begin{align*}
  f=\sum_{\uw\in \WSlenw} c_{\uw} f_{\uw}
\end{align*}
such that $\pi_w(f)=0$. 
Choose $w_1\in W^S$ minimal such that $c_{\uw}\neq 0$ for some
$\uw=(w_1,\dots,w_n)$. The relation $\pi_w(f)=0$ implies
$\pi_{w_1}(f)=0$. By induction hypothesis we may assume $w_1=w$. Hence
in $\sqgp/Q_w$ one again obtains the relation (\ref{nullteiler}),
where not all coefficients $c_{\uw'}$ vanish. As $\sqgp/Q_w$ is an
integral domain one obtains a contradiction to 1.~and 3.~for $n{-}1$.
Hence the standard monomials $\{f_{\uw}\,|\,\uw\in \WSlenw\}$
are linearly independent in $\pi_w(V(n\lambda)^\ast)$. As $Q_w$ is an
ideal all $f_{\uw}$ such that $\uw\in \WSlen\setminus \WSlenw$ belong
to $Q_w$. Hence by 2.~the set  $\{f_{\uw}\,|\,\uw\in \WSlenw\}$
is a basis of $\pi_w(V(n\lambda)^\ast)$. 

To verify 1.~note that $\Iqwld\subset Q_w$ because $Q_w$ is an
ideal. On the other hand $\dim(Q_w\cap V(n\lambda)^\ast)\le
\dim(\Iqwld\cap V(n\lambda)^\ast)$ by 3.
$\blacksquare$

\medskip

\begin{remark} \upshape
 Statement 1.~of the above proposition reproduces \cite[Th\'eor\`eme
 3]{a-Joseph95} for minuscule $\lambda\in \Ppluspi$. Note that Joseph's
 proof relies on a specialisation argument and hence yields the above
 result only for transcendental $q$. However, assuming here $\lambda$
 to be minuscule we only consider a very special case. In particular, it
 doesn't generally hold that the standard monomials contain a basis
 of $\ker(\pi_w)$.
 
 The second part of the above proposition reproduces
 \cite[3.8]{a-LaksResh92} for minuscule $\lambda$, again avoiding specialisation
 arguments. Finally, the third part of the above proposition is claimed
 in \cite[4.7]{a-LaksResh92} again for arbitrary $\lambda$ but only
 for transcendental $q$. $\square$
\end{remark}

\section{The AS-Cohen-Macaulay property}\label{ASCM-property}
The AS-Cohen-Macaulay property appeared in noncommutative algebraic
geometry as a graded analogue of the Cohen-Macaulay property for
commutative local
rings \cite{a-vdB97}, \cite{a-Jorg99}. As we are mainly interested in quantum flag manifolds we will
formulate this notion not in full generality but only for a slightly
restricted class of graded algebras. We refer the reader to
\cite{a-JorgZhang00} for the general definition and more details.

For any algebra $B$ we write $B^\circ$ to denote the opposite
algebra, i.e.~the vector space $B$ with multiplication defined by
$m(a,b)=ba$ for all $a,b\in B^\circ$.

\subsection{Gelfand-Kirillov dimension and depth}
Throughout this section we assume that $A=\bigoplus_{k=0}^\infty A_k$ is a
noetherian, $\N_0$-graded, connected algebra over a field $\field$ and
that $A$ is generated by the finite dimensional subspace $A_1$.
Define $V_n=\bigoplus_{k=0}^n A_n$ and let
\begin{align*} 
  \gkdim A =\limsup \frac{\log \dim_\field(V_
n)}{\log n}
\end{align*}
denote the Gelfand-Kirillov dimension of $A$.
\begin{example}\label{GKflags}\upshape
  For quantised flag manifolds one has $\gkdim(\sqgp)=l(w^S)+1$ where
  $w^S$ is the longest element in $W^S$. This follows from the
  corresponding result for $q=1$ because for any $\mu\in \Ppluspi$ the $\uqg$-module $V(\mu)$ has the
  same dimension as the corresponding $U(\gfrak)$-module. For the
  classical result consult e.g.~\cite[Exercise 2.6]{b-Ha} and recall that
  for finitely generated commutative $\field$-algebras the GK-dimension coincides with the Krull
  dimension \cite[Theorem 4.5 (a)]{b-KrauLen00}. $\square$
\end{example}
For any finitely generated left $A$-module $M$ the depth of $M$ is
defined by
\begin{align*}
  \depth_A M=\inf \{i\in \N\,|\, \Ext^i_A(\field,M)\neq 0\} \in \N\cup\{\infty\}.  
\end{align*}
Throughout this section we assume that there exists a sequence
$x_1,\dots,x_r$ of elements in $A_1$ such that the image of $x_i$ in
  $A/\langle x_1,\dots,x_{i-1}\rangle$ is normal for all
$i=2,\dots,r$ and such that $A/\langle x_1,\dots,x_r\rangle$ is a finite dimensional
$\field$ vector space.
As noted in \cite[p.~392]{a-Zhang97} the existence of such a normalising
sequence implies that $A$ has enough normal elements in the sense of
\cite[p.~392]{a-Zhang97}, \cite[Definition 2.13]{a-LenRi06}.
Under the previous assumptions we may use the following
definitions.
\begin{definition}\label{CMc-defi}
  Let $A$ be a noetherian, $\N_0$-graded, connected $\field$-algebra
  which has enough normal elements.
  \begin{enumerate}
    \item The algebra $A$ is called  AS-Cohen-Macaulay if {\upshape
      \begin{align}\label{CMc-eq}
        \depth_A A=\gkdim{A}=\depth_{A^\circ}A^\circ.
      \end{align}}
    \item The algebra $A$ is called AS-Gorenstein if it is of finite
    left and right injective dimension.
  \end{enumerate}
\end{definition}
The equivalence of the first definition with the standard definition of the
AS-Cohen-Macaulay property in terms of local cohomology follows for
instance from \cite[Remark 2.2.1]{a-LenRi06}. Similarly, the
equivalence of the above notion of AS-Gorensteinness with
\cite[Definition 0.2]{a-JorgZhang00} holds by \cite[Proposition 2.3(2)]{a-Zhang97}.

\begin{remark}\upshape
  Let $x_1,\dots,x_{\dim V(\lambda)}$ be a weight basis of
  $V(\lambda)^\ast$ such that $\wght(x_i)<\wght(x_j)$ implies
  $i<j$. It follows from the commutation relation (\ref{flag-rels})
  and (\ref{R-mat-prop}) that the
  image of this sequence in $\sqwgp$ yields a normalising sequence such
  that  $\sqwgp/\langle x_1,\dots,x_{\dim V(\lambda)}\rangle\cong
  \field$. Moreover, by Lemma \ref{noetherlem} the algebra $\sqwgp$
  is noetherian. Hence $\sqwgp$ satisfies all the assumptions in
  Definition \ref{CMc-defi} and it remains to verify (\ref{CMc-eq}) in
  order to obtain the AS-Cohen-Macaulay property. $\square$
\end{remark}

\subsection{Quantum algebras with a straightening law}
We recall the following definitions from \cite[1.1]{a-LenRi06}. The
second definition is a quantum version of a well known notion from
commutative algebra \cite[4.A]{b-BV88}, introduced and developed by C.~De Concini,
D.~Eisenbud, and C.Procesi. Consult \cite{a-DEP82}, \cite[4.E, 5.F]{b-BV88} for references to
the original literatur.
\begin{definition}
Let $A$ be an $\N_0$-graded algebra over a field $\field$ and $\Pi$ a finite
subset of $A$ equipped with a partial order $\le$. A standard monomial
on $\Pi$ is an element of $A$ which is either $1$ or of the form
$f_1\dots f_s$, for some $s\ge 1$, where $f_1,\dots,f_s\in \Pi$ and
$f_1\le\dots\le f_s$.
\end{definition}
\begin{definition}\label{ASL}
Let $A$ be an $\N_0$-graded algebra over a field $\field$ and $\Pi$ a finite subset of
$A$ equipped with a partial order $\le$. The algebra  $A$ is called a
quantum graded algebra with a straightening law on the poset $(\Pi,\le)$
if the following conditions are satisfied.
\begin{enumerate}
\item The elements of $\Pi$ are homogeneous with positive degree.
\item The elements of $\Pi$ generate $A$ as a $\field$-algebra.
\item The set of standard monomials on $\Pi$ is linearly independent.
\item If $f,g\in \Pi$ are not comparable for $\le$, then $fg$ is a
  linear combination of terms $F$ or $FG$, where $F,G\in \Pi$, $F\le
  G$, and $F<f,g$.
\item For all $f,g\in \Pi$, there exists $c_{fg}\in \field\setminus\{0\}$ such
  that $fg-c_{fg}gf$ is a linear combination of terms $F$ or $FG$,
  where $F,G\in \Pi$, $F\le G$, and $F<f,g$.
\end{enumerate}
To shorten notation we also call $A$ a quantum graded ASL on $\Pi$.
\end{definition}
The relevance of the above notions for our purposes stems from the
following quantum version of \cite[5.14]{b-BV88}.
\begin{theorem}\label{ASCM}{\upshape \cite[Theorem 2.2.5]{a-LenRi06}} If $A$ is a
  quantum graded ASL on a wonderful poset, then $A$ is AS-Cohen-Macaulay.
\end{theorem}

\subsection{The main result}
We are now in a position to prove the desired generalisation of
\cite[Theorem 3.4.4]{a-LenRi06}. 
Recall that $W^S_{\le w}=\{w'\in W^S\,|\,w'\le w\}$. 
\begin{proposition}\label{qgASL}
 Let $\lambda=\omega_s$ be a minuscule weight and
 $S=\pi\setminus\{\alpha_s\}$. For any $w\in W^S$ the algebra $\sqwgp$
 is a quantum graded ASL on the poset
 $\WSlew$ with the inverse Bruhat order.
 In particular $\sqgp$ is a quantum graded ASL on $W^S$.
\end{proposition}
\noindent {\bf Proof:} The generators of the algebra $\sqwgp$ can be identified
with the poset $\WSlew$. They are homogeneous of degree one and by
Proposition \ref{standardMinProp} the standard monomials are a
linearly independent set. Property 4 of Definition \ref{ASL} is
obtained analogously to the proof of \cite[Theorem 3.3.8]{a-LenRi06}.
More explicitly, let $w',w''\in \WSlew$ be not comparable in the
Bruhat order. By Proposition \ref{standardMinProp} one can write in
$\sqwgp$
\begin{align}\label{aiiAusdruck}
  f_{w'}f_{w''}=\sum_{j=1}^n a_jf_{w_j}f_{w_j'}
\end{align}
for some $w_j, w_j'\in \WSlew$, such that $w\ge w_j\ge w_j'$ for all
$j$, and $a_j\in \field\setminus \{0\}$. For any $j=1,\dots, n$ apply
the projection $\sqwgp\rightarrow S_q^{w_j}[G/P_S]$ to both sides of
(\ref{aiiAusdruck}). By Proposition \ref{standardMinProp} the image of
the right hand side under this projection is nonzero. Hence so is the
image of the left hand side. This implies $w'\le w_j$ and $w''\le w_j$
for all $j=1,\dots,n$. This proves property 4 with respect to the
inverse Bruhat order.

It remains to prove property 5. Let $w',w''\in \WSlew$. If $w'$ and
$w''$ are not comparable then the previous step shows that property 5
holds for $f_{w'}$ and $f_{w''}$. If on the other hand $w'>w''$ then
$w'\omega_s<w''\omega_s$. Hence by relations (\ref{flag-rels}) and
(\ref{R-mat-prop-dual}) as well as Proposition \ref{wonderprop}.4
there exist $w_j',w_j''\in W^S$, $j=1,\dots,n$, and $a_j\in \field$
such that
$w_j'>w'$, $w_j''<w''$ and 
\begin{align*}
  q^{(\omega_s,\omega_s)}f_{w''}f_{w'}=q^{(w'\omega_s,w''\omega_s)}f_{w'}f_{w''}+\sum_{j=1}^n 
  a_jf_{w_j'}f_{w_j''}.
\end{align*}
This concludes the proof of property 5 and hence the proof of the proposition. 
$\blacksquare$

\medskip

Proposition \ref{qgASL} implies an analogous result for quantum big
cells and determinantal varieties. To this end define $W^{S,e}=\{w'\in
W^S\,|\,e<w'\}$ and $W^{S,e}_{\le  w}=\{w'\in W^S\,|\,e<w'\le w\}$.
As the highest weight vector $f_e\in V(\lambda)^\ast$ $q$-commutes
with all generators of $\sqgp$ one immediately obtains the following result.
\begin{cor}
 Let $\lambda=\omega_s$ be a minuscule weight and
 $S=\pi\setminus\{\alpha_s\}$. For any $w\in W^S$ the algebra $\sqwgp(f_e^{-1})_0$
 is a quantum graded ASL on the poset
 $W^{S,e}_{\le w}$ with the inverse Bruhat order.
 In particular $\sqgp(f_e^{-1})_0$ is a quantum graded ASL on $W^{S,e}$.
\end{cor}

Using Theorem \ref{ASCM} and Remark \ref{invBruRemark} one now obtains the AS-Cohen-Macaulay
property for quantum flag manifolds of minuscule weight, for their
big cell, their Schubert varieties, and for their determinantal varieties.
\begin{cor}\label{minuscule-ASCM}
 Let $\lambda=\omega_s$ be a minuscule weight and
 $S=\pi\setminus\{\alpha_s\}$. For any $w\in W^S$ the graded algebras
 $\sqwgp$ and $\sqwgp(f_e^{-1})_0$
 are AS-Cohen-Macaulay. In particular $\sqgp$ and $\sqgp(f_e^{-1})_0$
 are AS-Cohen-Macaulay.
\end{cor}

\subsection{The AS-Gorenstein property}\label{ASG}
  Let $A$ be a noetherian, $\N_0$-graded, connected, AS-Cohen-Macaulay algebra with enough
  normal elements. In this case, by Stanley's Theorem
  \cite[Theorem 6.2]{a-JorgZhang00}, the algebra $A$ is AS-Gorenstein if and
  only if there exists $m\in \N_0$ such that the Hilbert series $H_A(t)=\sum_n \dim(A_n)t^n$ of $A$
  satisfies the functional equation 
  \begin{align}\label{pali}
    H_A(t)=\pm t^{-m}H_A(t^{-1})
  \end{align}
  as a rational function over $\Q$. For the quantum
  algebras considered in this paper the Hilbert series coincide with
  the Hilbert series of the corresponding commutative algebras. Hence
  these algebras are AS-Gorenstein if and only if their commutative
  analogues are Gorenstein (cf. \cite[Remark 2.1.10(ii)]{a-LenRi06}).
  We use Stanley's Theorem to verify the
  AS-Gorenstein property for $\sqgp$ and for its big cell
  $\sqgp(f_e^{-1})_0$ in the minuscule case.
  \begin{proposition}\label{cell-ASGor}
     Let $\lambda=\omega_s$ be a minuscule weight and
     $S=\pi\setminus\{\alpha_s\}$. Then the graded
     algebra $\sqgp(f_e^{-1})_0$ is AS-Gorenstein.
  \end{proposition}
  \noindent{\bf Proof:}
    It follows from Proposition \ref{dualprop}, the PBW-theorem for
    $\uqg$, and \cite[Proposition 4.2]{a-Kebe99} that the Hilbert
    series of $\sqgp(f_e^{-1})_0$
    coincides with the Hilbert series of a commutative polynomial ring
    generated by finitely many elements in positive degrees. Hence the
    Hilbert series of $\sqgp(f_e^{-1})_0$ satisfies condition (\ref{pali}).
   $\blacksquare$

  \medskip
  As we weren't able to locate a proof of the Gorenstein property of
  the commutative analogue of $\sqgp$ in the literature we give a
  proof for the minuscule case, following  methods of Stanley's
  paper \cite[Corollary 4.7]{a-Stanley78}.
  \begin{proposition}
     Let $\lambda=\omega_s$ be a minuscule weight and
     $S=\pi\setminus\{\alpha_s\}$. Then the graded
     algebra $\sqgp$ is AS-Gorenstein.
  \end{proposition}
  \noindent{\bf Proof:}
    It follows from Proposition \ref{charProp} that the Hilbert series
    of $\sqgp$ for $(\gfrak=\bfrak_n, \omega_s=\omega_n)$ coincides
    with the Hilbert series for $(\gfrak=\dfrak_{n+1},
    \omega_s=\omega_{n+1})$. Similarly, for  $(\gfrak=\cfrak_n,
    \omega_s=\omega_1)$ the Hilbert series of $\sqgp$ coincides with
    the Hilbert series of a polynomial ring. Hence we may restrict to
    the case where $\gfrak$ is of type ADE. 

    Let $\rho$ denote the half sum of all positive roots and define
  $\Delta^+_S=\{\alpha\in \Delta^+\,|\,(\omega_s,\alpha)=0\}$. It
  follows from the Weyl character formula that the Hilbert series of
  $\sqgp$ is given by
  \begin{align}\label{p(n)}
    H(t)=\sum_{n=0}^\infty p(n)t^n =p\big(t\frac{d}{dt}\big)\frac{1}{1-t}
  \end{align} 
  where $p$ denotes the polynomial defined by
  \begin{align*}
    p(n)= \prod_{\alpha\in \Delta^+\setminus \Delta^+_S}\frac{n+(\rho,\alpha)}{(\rho,\alpha)}. 
  \end{align*}
  It is straightforward to check that the function $H$ as a rational
  function of $t$ satisfies the relation
  \begin{align}\label{p(-n)}
    H(t^{-1})=p\big(-t\frac{d}{dt}\big)\frac{1}{1-t^{-1}}=-\sum_{n=1}^\infty p(-n)t^n.
  \end{align}
  Define $r=(\rho, \alpha_0)$ where $\alpha_0\in \Delta^+$ denotes the
  longest root, and for $k=1,\dots,r$ set $c_k=\#\{\alpha\in
  \Delta^+_S\,|\,(\rho,\alpha)=k\}$. Note that $p(-k)=0$ for all
  $k=1,\dots,r$. Moreover, $c_k=c_{r+1-k}$ and hence
  $p(n)=(-1)^{l(w^S)}p(-n-r-1)$ holds for all $n\in \Z$. In view of 
  (\ref{p(n)}) and (\ref{p(-n)}) this implies
  \begin{align*}
    H(t^{-1})=(-1)^{l(w^S)+1}t^{r+1} H(t)
  \end{align*}
  and the AS-Gorenstein property follows again from Stanley's Theorem.
  $\blacksquare$

\subsection{An example}\label{example}
  We end this paper with a first glance at what might
  happen in the the non-minuscule case. We verify the
  AS-Cohen-Macaulay property for $\gfrak=\bfrak_n=\sofrak_{2n+1}$ and
  $\lambda=\omega_1$ the first fundamental weight in the conventions
  of \cite[p.~58]{b-Humphreys}. In this case 
\begin{align*}
  V(\omega_1)\ot V(\omega_1)= V(2\omega_1)\oplus
  V(b_n\omega_2)\oplus V(0)
\end{align*}
where $b_n=2$ if $n=2$ and $b_n=1$ if $n\ge 3$. The quotient
of the tensor algebra $T(V(\omega_1)^\ast)$ by the ideal generated by
the subspace
$V(b_n \omega_2)^\ast$ of $ T^2(V(\omega_1)^\ast)$ coincides with
the well known quantum Euclidean space $O_q^{2n+1}(\field)$ introduced
in \cite{a-FadResTak1} for $\field=\C$ (cp.~also
\cite[9.3.2]{b-KS}). The algebra $O_q^{2n+1}(\field)$ contains a
central element $z$ which spans the component $V(0)^\ast \subset
T^2(V(\omega_1)^\ast)$. By construction one has
\begin{align}\label{quotrel}
  \sqgp\cong O_q^{2n+1}(\field)/\langle z \rangle.
\end{align}
It follows from the explicit relations \cite[Definition
  12]{a-FadResTak1}, \cite[9.3.2]{b-KS} that $O_q^{2n+1}(\field)$ is a
  quantum graded algebra with a straightening law on the poset
  $\{1,2,\dots,2n+1\}$. Hence Theorem \ref{ASCM} implies
\begin{align*}
  \depth_{O_q^{2n+1}(\field)} O_q^{2n+1}(\field)=\gkdim  O_q^{2n+1}(\field)
  =\depth_{O_q^{2n+1}(\field)^\circ} O_q^{2n+1}(\field)^\circ. 
\end{align*}
Note that $\gkdim O_q^{2n+1}(\field)=2n+1$ and $\gkdim \sqgp=2n$ by Example \ref{GKflags}.
In view of (\ref{quotrel}) it now follows from \cite[Lemma 2.1.7
  (iii)]{a-LenRi06} that $\sqgp$ is AS-Cohen-Macaulay. A similar
argument proves the AS-Cohen-Macaulay property for the corresponding
Schubert varieties. The quantum coordinate ring of the big cell is
isomorphic to $O_q^{2n-1}(\field)$ if $n\ge 3$. Hence
quantum big cells and determinantal varieties are also AS-Cohen-Macaulay in this case.
AS-Gorensteinness of $\sqgp$ and its big cell follows again directly
from Stanley's Theorem.

\medskip

\begin{center}{\bf Acknowledgement}
\end{center}

{\small The present work owes much of its existence to the article
\cite{a-LenRi06}. The author thanks T.H.~Lenagan and L.~Rigal
for detailed discussion and helpful comments.}


\providecommand{\bysame}{\leavevmode\hbox to3em{\hrulefill}\thinspace}
\providecommand{\MR}{\relax\ifhmode\unskip\space\fi MR }
\providecommand{\MRhref}[2]{%
  \href{http://www.ams.org/mathscinet-getitem?mr=#1}{#2}
}
\providecommand{\href}[2]{#2}

\textsc{Stefan Kolb, School of Mathematics and Maxwell Institute for
Mathematical Sciences, The University of Edinburgh, JCMB, The King's Buildings, Mayfield Road, Edinburgh
         EH9 3JZ, UK}
         
 \textit{E-mail address:} \texttt{stefan.kolb@ed.ac.uk}

\end{document}